\documentclass[a4paper,12pt,twoside]{article}
\usepackage[cp1251]{inputenc}
\usepackage[english,russian]{babel}
\inputencoding{cp1251}
\usepackage{amsmath}
\usepackage{graphicx}

\def\bb{\begin{equation}}
\def\ee{\end{equation}}

\def\sgn{\hbox{sgn}}

\def\l{\lambda}
\def\O{\Omega}
\def\o{\omega}

\numberwithin{equation}{section}

\begin{document}
%\large
%\title{ОТ СЛАБЫХ РАЗРЫВОВ К БЕЗДИССИПАТИВНЫМ УДАРНЫМ ВОЛНАМ}
%\author{Р.Н.Гарифуллин, Б.И.Сулейманов\footnote{ E-mail:rustem@matem.anrb.ru, bisul@mail.ru }}
%\date{\today}
\begin{center}
{\bf\large НЕКОТОРЫЕ ОСОБЕННОСТИ ИЗГИБОВ 
СТЕРЖНЯ ПРИ CИЛЬНОМ ПРОДОЛЬНОМ СЖАТИИ}
\\
%\vspace*{3mm}
{\it А.А.Ершов, Б.И.Сулейманов}  \footnote{ E-mail:ale10919@yandex.ru,
 bisul@mail.ru }
\\ {\it Институт механики и Математики УрО РАН(отдел динамических систем); Институт математики с ВЦ УНЦ РАН\\
$450008$, Уфа, Россия}\\
%\\ \vspace*{3mm}
\end{center}
%------------------------- A B S T R A C T --------
\begin{abstract} 
Изучаются некоторые типичные процессы изгибов стержня при сильном продольном сжатии. Соответствующее динамическое уравнение изгибов рассматривается как возмущение двумерного уравнения Лапласа. % с малым параметром. %при производной четвертого порядка по одной из переменных. 
 Установлено, что для данных процессов расширяющиеся области быстрых нарастаний изгибов начинаются в малых окрестностях точек сингулярностей решений предельного уравнения Лапласа. Начальные стадии этих нарастаний описаны с помощью интеграла Харди. 
%Подробно описано поведение этих эталонных интегралов. 

\end{abstract}

%\section{Введение}
\section*{\small{ ВВЕДЕНИЕ}}

{\bf 0.1.} Динамическое уравнение изгибов стержня 
 \cite[\S17] {Vol}, \cite[\S15] {PaG}, \cite{LaI}
 %(выведенном, согласно принципу Даламбера учета действия 
%инерционных сил)
\begin{equation} \label{uvst}
Pu''_{xx}+\rho u''_{tt}+Gu''''_{xxxx}=0\end{equation}
с постоянными $P$, $\rho$, $G$
($P$ -- величина продольной сжимающей силы, $\rho$ -- линейная плотность стержня,
 $G$ -- жесткость на изгиб стержня) после замены %перехода к временной переменной 
\begin{equation} \label{scalet}
y=(\frac{P}{\rho})^{1/2}t\end{equation}
принимает вид возмущения уравнения Лапласа
\begin{equation} \label{sust}
u''_{xx}+u''_{yy}+\varepsilon^2 u''''_{xxxx}=0\end{equation}
с параметром $\varepsilon=(G/P)^{1/2}$, который при $P>>G$ %(то есть при малой жесткости $G$, либо при больших значениях сжимающей силы  $P$)
 является малым. Посредством этого %же  
уравнения
моделируются и ряд других процессов
 \cite{Il}, \cite{Pa}, \cite[\S 30]{Vo},  \cite{Ah}.  

В отличие от %решений задачи Коши для 
некорректного по Адамару 
 предельного уравнения Лапласа  
\begin{equation} \label{lapl}
U''_{xx}+U''_{yy}=0,\end{equation}
%решение начальной задачи для 
это его возмущение корректно: %при фиксированных значениях $y$ и $\varepsilon$
 модуль решения (\ref{sust}) 
\begin{equation} \label{modi}
u=\exp(\lambda (1-\varepsilon^2 \lambda^2)^{1/2}y)\sin(\lambda x) 
\end{equation}
не раcтет при $\lambda \to \infty $, так как при переходе  $\lambda$ через значение
\begin{equation} \label{inkr}
\lambda=\frac{1}{\varepsilon}
\end{equation}
показатель экспоненты в (\ref{modi}) cтановится чисто мнимым. 
%(Это значение формулой $ P=G\lambda^{-2}$, где $\lambda=\lambda_*$ %соответствует низшей из гармоник  $\cos(\lambda x)$, задает критическую %нагрузку стержня -- при ее превышении возникает, как  явление {\it флаттера} %\cite[\S15] {PaG}, так как амплитуда этой гармоники экспоненциально растет %по времени.) 

{\bf 0.2.} Данное обстоятельство, заметим, иллюстрирует насыщающее  действие малой дисперсии на неустойчивость системы приближения нелинейной геометрической оптики 
%(НГО)  
\cite{s46} 
\begin{equation} \label{obosg}
h'_{\tau}+(h v)'_x=0,\quad v'_{\tau}+vv'_x-\alpha(h)h'_x=0 \end{equation}
в %его  
применении к известному нелинейному уравнению Шредингера
%(НУШ)
\begin{equation} \label{Schr}-i\mu F'_{\tau}=\mu^2F''_{xx}+K(|F|^2)F, \qquad \alpha(h)=2K'(h)>0, \qquad 0<\mu<<1. \end{equation}
%$\alpha(\h)=2K'(\h)>0$.

В самом деле, %данное приближение, описываемое эллиптической системой (\ref{obosg})
%\begin{equation} \label{obosg}
%\h'_{\tau}+(\h v)'_x=0,\quad v'_{\tau}+vv'_x-\alpha(\h)\h'_x=0, \end{equation}
%($h$---плотность, $v$---скорость, $\alpha(h)>0$),
%в применении к (\ref{Schr}) возникает  cледующим образом: 
подстановка в  (\ref{Schr}) выражения
$F = h ^{1/2} \exp(i\varphi/\mu)$ 
%\end{equation}
и дифференцирование по переменной $x$ второго уравнения результата этой подстановки
$$h '_{\tau}+ 2(h \varphi'_x)'_ x = 0,\quad    \varphi'_{\tau} + (\varphi'_x)^2 - K(h ) = \mu^2\frac{( \sqrt{h })''_{xx}}{\sqrt h }   $$
дает систему $(v= 2(\varphi)'_x)$
\begin{equation} \label{geomg}
h '_{\tau} + (\ v)'_ x = 0,\quad    v'_{\tau} + vv'_x-\alpha(h )h '_x =
2\mu^2(\frac{( \sqrt{h })''_{xx}}{\sqrt h })'_x,
\end{equation}
которая при 
$\mu = 0$ cовпадает с неустойчивой системой уравнений (\ref{obosg}). 
Рассмотрим начальную стадию эволюции решений системы (\ref{geomg}), когда возмущения еще малы, считая, что невозмущенному состоянию среды отвечает частное решение этой системы $v=0$, $h=1$, и пусть $\hat{v}$, $\hat{h}$ --- малые возмущения соответствующих компонент решения. Полагая в (\ref{geomg}) $v=\hat{v}$, $h=1+\hat{h}$ и отбрасывая старшие по степеням этого возмущения члены, получим для поправок линейную систему %($y=(\alpha(1))^{1/2} \tau$ , $\varepsilon=(\alpha(1))^{-1/2} \mu$)
\begin{equation} \label{geogd}
\hat{h'}_{\tau}+ \hat {v'}_ x = 0,\quad    \hat{v'}_{\tau} -\alpha(1){\hat {h} }'_x =
\mu^2{\hat{h} '''}_{xxx},
\end{equation}
компоненты решения которой после перехода к независимой переменной  \linebreak $y=(\alpha(1))^{1/2} \tau$  задают решения $\hat{h}$, $\hat{v}$ уравнения (\ref{sust}) c параметром 
$\varepsilon=\alpha(1)^{-1/2} \mu$. Тогда как решение начальной задачи для предельной к (\ref{geogd}) cистемы 
%\begin{equation} \label{geopd}
$$
\hat{h'}_{\tau}+ \hat {v'}_ x = 0, \quad \hat{v'}_{\tau} -\alpha(1 ){\hat{h} }'_x =0,
$$
%\end{equation}
некорректно по Адамару (это обстоятельство демонстрирует \cite[\S 3] {ZhT} неустойчивость системы приближения нелинейной геометрической оптики (\ref{obosg})), гиперболическая \cite [Гл. III, \S 3.1] {Gel'fSh} система (\ref{geogd}) лишена этого недостатка.

{\bf 0.3.} Значение (\ref{inkr}) формулой $P_E=G\lambda_*^{2}$, где $\lambda_*$ соответствует низшей из возможных гармоник, определяет классическую критическую нагрузку Эйлера. При ее превышении амплитуда низшей из гармоник (\ref{modi}) экспоненциально растет по времени. Поэтому иногда высказывается мнение о том, что изучение решений (\ref{uvst}) при значениях $P$, превышающих критическое, с практической точки зрения неактуально: `` в конце концов все равно, произойдет ли ставшим неизбежным разрушение стержня через одну или две секунды... '' \cite[\S15] {PaG}. 

Необоснованность этого мнения следует хотя бы из следующего комментария Ишлинского к статье \cite{LaI}, в которой среди членов ряда Фурье, задающего решения уравнения (\ref{uvst}), в первую очередь рассматривались как раз экспоненциально растущие по времени: 
%``Превращение продольной трубы в гофрированную 
``...представлялось тогда загадкой ...%, ровно как и 
появление высших форм потери устойчивости при внезапном (ударном) сжатии первоначально прямолинейного стержня. Следует пояснить, что начало потери устойчивости возникает еще при упругих деформациях. В дальнейшем, как правило, развивается пластичность и именно остаточные деформации наблюдаются по окончании...
%процесса обжатия трубы или 
продольного сжатия стержня'' \cite [c.5] {Isch}. 

Вывод \cite{LaI} о том, что среди членов рядов Фурье, представляющих решения уравнений (\ref{sust}), наиболее быстро растут те, что соответствуют гармоникам (\ref{modi}) со значениями $\lambda$, близкими к точке максимума 
$1/(\sqrt{2}\varepsilon)$
 функции $\lambda(1-\varepsilon^2\lambda^2)^{1/2},$
%(то есть, с номерами, которые примерно равны $2/3$ от номера наивысшей возможной формы потери устойчивости, определяемого критическим значением  (\ref {inkr}))
 неоспорим: со временем амлитуды %$R(\lambda)\exp(\lambda(1-\varepsilon^2\lambda^2)^{1/2})$ 
 этих членов рядов Фурье становятся значительно больше других. Он подтвержден экспериментами по сжатию стержней как посредством подрывов у одного из их концов, так и при ударах молотком \cite[ c.547] {Ishl}. Результаты проведенного нами численное моделирование решений конкретных задач, относящихся к кругу тех, что рассматриваются в этой статье, также согласуются с этим выводом. 

Интерес, однако, представляют и другие аспекты поведения решений (\ref{sust}), исследованию которых посвящен уже целый ряд работ (помимо процитированных выше, укажем еще на прешественницу \cite{Kon} cтатьи \cite{LaI}, а также на публикации последних лет \cite{Mai} -- \cite{Pai} -- в части из них речь идет о ранних этапах эволюции решений уравнения (\ref{sust})). Мы же здесь хотим обратить внимание на следующую сторону проблемы описания динамики сжимаемого стержня:

при большой приложенной силе изгибание стержня и вступление затем в силу эффектов пластичности происходит очень быстро. Поэтому, на наш взгляд, в такой ситуации интересоваться надо именно 
%не 
асимптотиками решений уравнения (\ref{uvst}) %при больших временах $t$, а их асимптотиками 
по малому параметру $\varepsilon$.
 Кстати, согласно (\ref{scalet}) при больших значениях $P$ малым -- порядка $O(\varepsilon)$ -- временам $t$ соответствуют конечные значения $y$.

{\bf 0.4.} Но до сих пор малость параметра $\varepsilon$ для описания поведения решений уравнения (\ref{uvst}) и эквивалентного ему уравнения (\ref{sust}) почти не использовалась 
(исключением является недавняя статья \cite{Pai}). Не рассматривались, в частности, стандартные ряды теории возмущений %по малому параметру $\varepsilon$
\begin{equation}\label{vnesh}u(x,y,\varepsilon)=U_0(x,y)+\varepsilon^2U_1(x,y)+\sum_{n=2}^{\infty}\varepsilon^{2n}U_n(x,y),\end{equation}  
главные члены $U_0(x,y)$ которых есть решения уравнения Лапласа (\ref{lapl}), а  остальные находятся из рекуррентной последовательности уравнений 
\begin{equation} \label{lapk}
(U_n)''_{xx}+(U_n)'_{yy}=-(U_{n-1})''''_{xxxx}.\end{equation}
%За исключением, быть может \cite[\S 3] {ZhT}, 
%Для описания решений (\ref{sust}) 
%Для данной цели 
%Пока почти не привлекались и просто решения предельного уравнения Лапласа (\ref{lapl}) ( исключение: недавняя работа \cite{Pai}). 
Это может объясняться следующими причинами:

1) само уравнение изгибов стержня (\ref{sust}) интегрируемо -- оно, в частности, допускает стандартную процедуру разделения переменных. И, как уже отмечено, % в отличие от предельного уравнения Лапласа,
 при $\varepsilon>0$ начальная задача для гиперболического уравнения (\ref{sust}) корректна. Впрочем, некорректность начальных задач для предельного уравнения Лапласа -- не такое уж неодолимое препятствие для их изучения. Такие задачи часто встречаются в приложениях \cite[\S3, пример 4] {Tih} и при наличии дополнительной информации об их решениях%(например, об ограниченности этих решений наперед заданной постоянной в некоторой области, большей, чем рассматриваемая), 
, в принципе, поддаются исследованию \cite{Tih} -- \cite{Leo}; 

2) решениям начальных задач для уравнения Лапласа (\ref{lapl}) при конечных $y$ присуща взрывная потеря пригодности даже тогда, когда их начальные данные 
\begin{equation} U(x,y_0)=A(x), \qquad U'_y(x,y_0)=B(x), \label{Kof}\end{equation}
есть аналитические функции \cite{Tarh}. При приближении к точкам сингулярностей решений задачи (\ref{lapl}), (\ref{Kof}) теряют пригодность и все члены рядов теории возмущений (\ref{vnesh}), определяемые данными решениями (см. \S1 настоящей работы). 

{\bf 0.5.} Однако, надо заметить, что несмотря на сингулярности, присущие \cite{s46}, \cite{ZhT},
\cite{Gursh} -- \cite{Dub} решениям %эллиптической
 системы уравнений  (\ref{obosg}), для изучения свойств решений уравнения (\ref{Schr}) данная неустойчивая система применяется давно и успешно. %, так как другие аналитические методы изучения их решений общего положения %при конечных значениях $x$ и $\tau$ 
%в настоящее время не просматриваются (за исключением разве что интегрируемого методом обратной задачи рассеяния случая $K(|G|^2)=|G|^2$). 
При этом типичные сингулярности решений (\ref{obosg}) играют важную роль и для описания характерных особенностей поведения решений уравнения (\ref{Schr}) c малой дисперсией: по
  крайней мере для 
%давно считается совершенно естественным и практикуется очень широко. 
двух типичных сингулярностей решений (\ref{obosg}) в настоящее время проведено весьма подробное исследование поправочного влияния малой дисперсии (правой части системы (\ref{geomg})) на эти две сингулярности: в случае  первой из этих сингулярностей, отвечающей процессу провального самообострения амплитуды $h$, такому исследованию были посвящены работы \cite{Kudb}, \cite{Bkud};  в случае второй сингулярности, отвечающей особенности квадратного корня комплексно -- аналитических отображений, это влияние для широкого ряда случаев описывается в довольно большом числе публикаций, среди которых ключевыми, как мы считаем, являются \cite{Dub}, \cite{Tovb}.

В данной работе для двух характерных сингулярностей решений уравнения Лапласа (\ref{lapl}) проводится исследование поправочного влияния на них члена $\varepsilon^2u_{xxxx}$ в уравнении (\ref{sust}). И тем самым дается %детальное 
описание начала некоторых типичных процессов резкого изменения в поведении решений (\ref{sust}), которое при отказе от использования рядов возмущений (\ref{vnesh}) оставалось ``за кадром''. 

Первая характерная сингулярность решений уравнения Лапласа, рассматриваемая нами в этой работе, задается линейными комбинациями вещественной и мнимой частей комплексно-аналитических функций с полюсами вида 
\begin{equation}\label{res}U=\frac{1}{(-y+ix)}+\dots \qquad(%при 
x^2+y^2\to 0).\end{equation}
(Жданов и Трубников \cite[\S 3] {ZhT} особо отмечают универсальность соответсвующих сингулярных решений, которые, как стоит отметить, в главном порядке выражаются через производные $\Phi'_x(x,y)$ и $\Phi'_y(x,y)$ {\it фундаментального} решения $\Phi(x,y)=\ln(x^2+y^2)^{1/2}/(2\pi)$ уравнения Лапласа.)
%$E''_{xx}+E''_{yy}=\delta(x,y).$ Важность и универсальность таких сингулярностей особо отмечается Ждановым и Трубниковым \cite[\S 3] {ZhT}) 

Вторая сингулярность в главном порядке есть вещественная %частью %, задающей %вещественной и 
%мнимую часть 
%частями 
часть корня 
\begin{equation}\label{root}U={(2(-y+ix))}^{1/2}.\end{equation}
(Данная особенность решения уравнения Лапласа типична с точки зрения математической теории катастроф \cite [1.8]{Gus}.) Возникающую задачу об истинном поведении решений полного уравнения (\ref{sust}) в окрестности точки $(x=0, y=0)$, которое при $y \leq 0$ вне малой по $\varepsilon$ окрестности этой точки задается соответствующим рядом (\ref{vnesh}), в известной мере можно рассматривать как линеаризованный вариант проблемы, исследовавшейся в работах \cite{Dub}, \cite{Tovb}.
%Как показано ниже в \S2, для таких сингулярностей перестройка соответствующих решений (\ref{sust}) в окрестностях начала координат описывается с помощью 
% специальных функций (
%\begin{equation}I_{k/2}(X,Y)=\int\limits_0^{\infty}\lambda^{k/2}\exp{(\lambda (Y-iX)-\lambda^3/3)} d\lambda, (\label{pOrlv}
%\end{equation}
%где растянутые переменные $X$ и $Y$ задаются формулами 
%\begin{equation}\label{skail} X=\frac{x}{\varepsilon^{2/3}}, \qquad %Y=\frac{y}{\varepsilon^{2/3}}.\end{equation}

В \S2 и \S3 показывается, что для двух этих характерных сингулярностей перестройка соответствующих решений (\ref{sust}) в окрестности начала координат описывается с помощью интеграла 
 {\it Харди} \cite{Hardy}, \cite [18.3.1, 21.1.3]{Rie}
\begin{equation}I_{\nu}(X,Y,m)=\int\limits_0^{\infty}k^{\nu}\exp{(k(Y-iX)-k^m/2)} dk ,\quad\mathrm{Re}\nu>-1, \label{pOrlv}
\end{equation}
где $m=3$, индекс $\nu$ принимает значения $-1/2, 0, 1/2, 3/2$,% а растянутые переменные $X$ и $Y$ задаются формулами 
\begin{equation}\label{skail} X=\frac{x}{\varepsilon^{2/3}}, \qquad Y=\frac{y}{\varepsilon^{2/3}}.\end{equation}
Из асимптотик $I_{\nu}(X,Y,3)$ при $X^2+Y^2\to \infty$ cогласно методу согласования \cite{Ilam} следует, что за точками рассматриваемых типичных сингулярностей решений уравнения Лапласа (\ref{lapl}) cоответствующие решения уравнения (\ref{sust}) экспоненциально растут при $\varepsilon\to 0$ -- именно в малых окрестностях этих точек начинаются расширяющиеся по мере роста $y$ области сильных изгибов стержня. 

{\bf 0.6.} В настоящей статье, как и в \cite{LaI}, %наряду с идеальным случаем начальной задачи для этого уравнения в данной работе 
исследуется случай свободно (шарнирно) опертого стержня, которому отвечают начально-краевые задачи на промежутке 
$-L\leq x\leq L$ $(L>0)$ с краевыми условиями 
\begin{equation} \label{kra}u(x,y)|_{x=\pm L}=u''_{xx}( x,y)|_{x=\pm L}=0.\end{equation}

\section{РЯД ТЕОРИИ ВОЗМУЩЕНИЙ}

{\bf1.1.} Для удобства сначала рассмотрим идеализированный случай чиcто начальной задачи для уравнений %изгибов стержня 
 (\ref{sust}) c начальными данными (\ref{Kof}), заданными в момент времени $y_0=-1$. В предположении, что решение  $U_0(x,y)$ уравненния  Лапласа (\ref{lapl}), удовлетворяющее начальным условиям (\ref{Kof}), существуют при $-2<y<0$, все 
коэффициенты $U_n(x,y)$ ряда (\ref{vnesh}), который представляет решение начальной задачи (\ref{Kof}),  явным образом выписываются через его главный член $U_0(x,y)$. Действительно, при $k>1$ эти коэффициенты есть решения рекуррентной системы уравнений (\ref{lapk}), удовлетворяющие начальным условиям
$$U_n(x,y)|_{y=-1}=(U_n(x,y))'_y|_{y=-1}=0.$$
%\begin{equation}\label{koshk}U_n(x,y)|_{y=-1}=(U_n(x,y))'_y|_{y=-1}=0. \end{equation}

Решая ее последовательно, находим, что 
\begin{equation} \label{rekk}
U_1(x,y)=\frac{1}{2}\left(\frac{(U_0(x,y))''_{yy}}{2}-\frac{(U_0(x,-y-2))''_{yy}}{2}-(y+1)(U_0(x,y))'''_{yyy}\right), $$

$$U_2(x,y)=\frac{1}{4}
\left(\frac{3(U_0(x,y))''''_{yyyy}}{4}-\frac{3(U_0(x,-y-2))''''_{yyyy}}{4}-\right.$$
$$\left.-(y+1)\left((U_0(x,y))'''''_{yyyyy}+\frac{(U_0(x,-y-2))'''''_{yyyyy}}{2}-\frac{y+1}{2}(U_0(x,y))''''''_{yyyyyy}\right)\right), $$
$$ U_n(x,y)=\sum\limits_{m=0}^{n}(y+1)^{m}a_m \frac{\partial^{(2n+m)}U_0(x,y)}{\partial y^{(2n+m)}}+\sum\limits_{m=0}^{n-1}(y+1)^{m}b_m\frac{\partial^{(2n+m)}U_0(x,-y-2)}{\partial y^{(2n+m)}},\end{equation}
где $a_m$ и  $b_m$ -- конкретные вещественные постоянные. Справедливость последней из формул (\ref{rekk}) при всех $n>2$ устанавливается по индукции.

В принципе, это же разложение можно получить, решая начальную задачу (\ref{sust}), (\ref{Kof}) c помощью преобразования Фурье, учитывая равенства 
$$\ch(k(-y-2+1))=\ch(k(y+1)), \qquad \sh(k(-y-2+1))=-\sh(k(y+1)), $$
и используя разложение функции $k(1-\varepsilon^2k^2)^{1/2}$ в ряд по степеням
% малого параметра 
 $\varepsilon$. 
 
Но как мы сейчас покажем, ряд (\ref{vnesh}), (\ref{rekk}) полезен не только для описаний решений чисто начальной задачи (\ref{sust}), (\ref{Kof}).

{\bf 1.2.} Существенной особенностью представлений (\ref{rekk}) коэффициентов $U_k(x,y)$
ряда (\ref{vnesh}) является то, что все они через его главный член $U_0(x,y)$ выражаются поcредством конечного числа дифференцирований {\it исключительно по переменной $y$} и умножений таких его производных на натуральные степени $(y+1)$.

Допустим теперь, что при $-2 \leq y<0$ функция $U_0(x,y)$ является гладким решением уравнения Лапласа (\ref{lapl}) для всех $-L\leq x \leq L$, которое на краях этого замкнутого промежутка обращается в нуль:
$$U_0(x,y)|_{x=\pm L}=0.$$
%\begin{equation} \label{krag}U_0(x,y)|_{x=\pm L}=0.\end{equation} 
Из удовлетворения функцией $U_0(x,y)$ уравнению Лапласа (\ref{lapl}) следует, что тогда на этих краях обращаются в нуль и вторые производные $(U_0(x,y))''_{xx}$:
$$(U_0(x,y))''_{xx}|_{x=\pm L}=0.$$
%\begin{equation} \label{krag}(U_0(x,y))''_{xx}|_{x=\pm L}=0.\end{equation} 
Поэтому и все остальные члены $U_n(x,y)$ ряда (\ref{vnesh}), определяемые формулами (\ref{rekk}), будут тогда автоматически удовлетворять краевым условиям
$$U_n(x,y)|_{x=\pm L}=(U_n(x,y))''_{xx}|_{x=\pm L}=0.$$
%\begin{equation} \label{krak}U_n(x,y)|_{x=\pm L}=(U_n(x,y))''_{xx}|_{x=\pm L}=0.\end{equation}
И, значит, в этих условиях при $-1\leq y<0$ и $-L\leq x \leq L$ ряд (\ref{vnesh}), (\ref{rekk}) представляет решение $U(x,y,\varepsilon)$ уравнения (\ref{sust}), удовлетворяющее как начальному условию (\ref{Kof}) (заданному в момент времени $y_0=-1$), так и {\it всем} условиям свободного опирания стержня (\ref{kra}).

Понятно, что все вышесказанное легко переносится на случай, например, начального момента времени $y_0=0$ и гладкости решения предельного уравнения на любом промежутке $-M\leq y<M$ $(M>0)$.

В дальнейшем для определенности мы будем говорить лишь о начальных данных, заданных при $y_0=-1$. Далее также всюду %, кроме Замечания ?,
 предполагается, что функция $U_0(x,y)$ есть решение уравнения Лапласа, которое является гладким на интервале $-2\leq y<0$ для всех $x$ (при рассмотрении чисто начальных задач) или для $-L\leq x \leq L$ (при рассмотрении начально-- краевых задач с краевыми условиями, заданными при $x=\pm L$).

{\bf 1.3.} Как уже говорилось во Введении, для решений  $U_0(x,y)$ начальных задач для уравнения Лапласа (\ref{lapl}) характерна взрывная потеря пригодности. Предположим, что она происходит в начале координат $(x=0, y=0)$ и описывается сингулярностями вроде, например, типичных сингулярностей, которые в главном при $x^2+y^2\to 0$ порядке опиcываются вещественными и мнимыми частями решений (\ref{res}), (\ref{root}). В частности, предполагается, что каждое дифференцирование по переменной $y$ ухудшает эту сингулярность на величину порядка $O((x^2+y^2)^{-1/2})$. Тогда из формул (\ref{rekk}) cледует, что по мере роста  номера $n$ коэффициенты $U_n(x,y)$ ряда теории возмущений ведут себя все хуже -- каждый из них ``сингулярнее'' предыдущего на величину порядка $O((x^2+y^2)^{-3/2})$. И ясно, что в малой окрестности начала координат ряд (\ref{vnesh}) теряет свою пригодность.

Согласно идеологии метода согласования \cite{Ilam} правильное поведение в этой окрестности соответствующих решений (\ref{sust}) описывается в растянутых переменных $X=\varepsilon^{-\alpha}x$, $Y=\varepsilon^{-\beta}y$, где положительные постоянные $\alpha$ и $\beta$ нужно подобрать так, чтобы  с одной стороны величины $\varepsilon^2(x^2+y^2)^{-3/2}$ при конечных значениях $X$ и $Y$ имели порядок $O(1)$, а с  другой, чтобы все члены предельного уравнения Лапласа (\ref{lapl}) после такого растяжения имели одинаковый порядок.

Таким образом, получаем, что правильный выбор растянутых переменных задается как раз формулами (\ref{skail}). После замен независимых переменных согласно этим формулам (\ref{sust}) переходит в уравнение 
$$u''_{XX}+u''_{YY}+\varepsilon^{2/3} u''''_{XXXX}=0,$$
%\begin{equation}\label{susi} u''_{XX}+u''_{YY}+\varepsilon^{2/3} %u''''_{XXXX}=0,\end{equation}
предельным к  которому по-прежнему является уравнение Лапласа
\begin{equation}\label{Lapl} u''_{XX}+u''_{YY}=0.\end{equation}

Поэтому для широкого ряда случаев, включающих рассматриваемые нами в этой работе случаи cингулярностей рядов (\ref{vnesh}), в главном по малому параметру $\varepsilon$ порядке перестройки соответствующих решений уравнения (\ref{sust}) будут задаваться решениями этого уравнения Лапласа (\ref{Lapl}) (ему, в частности, удовлетворяет интеграл Харди (\ref{pOrlv})). Только,  в отличие от сингулярных в начале координат решений $U_0(x,y)$ уравнения (\ref{lapl}), они будут глобально гладкими по переменным $X$ и  $Y$ решениями. Асимптотики их при $X^2-Y\to \infty$ определяются особенностями рядов (\ref{vnesh}), (\ref{rekk}) в начале координат. 

\section{ПОЛЮСНАЯ CИНГУЛЯРНОСТЬ}

{\bf 2.1.} Вещественные решения $U_0(x,y)$ уравнений Лапласа (\ref{lapl})
с типичными особенностями, которые в главном порядке задаются вещественной и мнимой частями решений вида (\ref{res}), при $x^2+y^2\to 0$ разлагаются в ряд
\begin{equation} \label{polm}U_0(x,y)
=\frac{-y}{y^2+x^2}+\sum\limits_{i,j=0}^{\infty}a_{ij}x^iy^j,\end{equation}
\begin{equation} \label{polmm}U_0(x,y)
=\frac{-x}{y^2+x^2}+\sum\limits_{i,j=0}^{\infty}b_{ij}x^iy^j.\end{equation}

Действуя согласно методу эталонных задач \cite{Bab}, рассмотрим следующие два частных случая таких решений 
%ИСПРАВЛЕНИЕ 1.
\begin{equation} \label{rskri}U_0(x,y)
\equiv\frac{-y}{{y^2+x^2}},\end{equation}
%А НЕ
%\begin{equation} \label{rskrn}U_0(x,y)
%\equiv\frac{-y}{{y^2+x^2}},\end{equation}
\begin{equation} \label{rskrr}U_0(x,y)
\equiv \frac{-x}{{y^2+x^2}}\end{equation}
и определяемые ими чисто начальные задачи для уравнения (\ref{sust}) с условиями 
\begin{equation} \label{koshn} u(x,-1)=A(x),\quad u'_y(x,-1)=B(x),\end{equation}
где 
\begin{equation} \label{koshAB} A(x)=\frac{x}{{1+x^2}}, \qquad B(x)=\frac{1-x^2}{(1+x^2)^2}
\end{equation} и
\begin{equation} \label{koshABr} A(x)=-\frac{1}{{1+x^2}}, \qquad B(x)=-\frac{2x}{(1+x^2)^2}.\end{equation}
Далее $u_{re}(x,y,\varepsilon)$ и $u_{im}(x,y,\varepsilon)$ будут обозначать решения эталонных начальных задачи (\ref{sust}), (\ref{koshn}) с начальными данными (\ref{koshAB}) и, соответственно, (\ref{koshABr}).

{\bf 2.2.} Эти эталонные задачи явно решаются при помощи преобразования Фурье, которое функции $v(x)$ cопоставляет интеграл
\begin{equation}\label{inpf}\widetilde{v}(k)=\int\limits_{-\infty}^{\infty}v(x)
\exp{(ikx)}dx.\end{equation}

После преобразования Фурье задача (\ref{sust}), (\ref{koshn}), (\ref{koshAB}) сводится \cite[Гл.8, \S 1.108, \S 1.124] {Bryp} к начальной задаче
\begin{equation}%\left\{
%\begin{array}{l}
\widetilde{u}_{yy}-(k^2-\varepsilon^2 k^4)\widetilde{u}=0,\quad -\infty<k<+\infty,\;y>-1,\label{oduf}\end{equation}
%\\
$$\widetilde{u}\big|_{y=-1}=\pi \exp{(-|k|)},\qquad
%\\
\widetilde{u}_y\big|_{y=-1}=\pi|k|\exp{(-|k|)},
%\end{array}\right 
$$
решение которой имеет вид 
%\begin{equation}\label{polfm}
$$\widetilde{u}(k,y,\varepsilon)=\pi \exp{(-|k|)}\ch\left(\sqrt{k^2-\varepsilon^2 k^4}(y+1)\right)
+\pi \exp{(-|k|)}\frac{\sh\left(\sqrt{k^2-\varepsilon^2 k^4}(y+1)\right)}{\sqrt{1-\varepsilon^2 k^2}}.$$
%\end{equation}
Применяя к этой функции обратное преобразование Фурье, получаем явное решение эталонной начальной задачи (\ref{sust}), (\ref{koshn}), (\ref{koshAB}) в виде 
\begin{equation}\label{polfm}
u_{re}(x,y,\varepsilon)=\frac{1}{2\pi}\int\limits_{-\infty}^\infty\widetilde{u}(k,y,\varepsilon)\exp{(-ikx)}dk
=\frac{1}{\pi}\mathrm{Re}\int\limits_0^\infty\widetilde{u}\exp{(-ikx)}dk=$$

$$=\mathrm{Re}\int\limits_0^{\infty}
\frac{\exp{(\sqrt{k^2-\varepsilon^2k^4}(y+1))}}{2}\left(1+\frac{1}{\sqrt{1-\varepsilon^2k^2}}\right)
\exp{(-k-ikx)}dk+$$
$$+\mathrm{Re}\int\limits_0^{\infty}
\frac{\exp{(-\sqrt{k^2-\varepsilon^2k^4}(y+1))
}}{2}\left(1-\frac{1}{\sqrt{1-\varepsilon^2k^2}}\right)
\exp{(-k-ikx)}dk.
\end{equation}

Преобразование Фурье начальных данных (\ref{koshABr}) также задается элементарными функциями 
\cite[Гл.8, \S 1.110, \S 1.122]{Bryp} -- образом Фурье решения задачи 
(\ref{sust}), (\ref{koshn}), (\ref{koshABr}) является решение уравнения (\ref{oduf}), которое удовлетворяет условиям 
$$\widetilde{u}\big|_{y=-1}=-i\pi \sgn{(k)}\exp{(-|k|)},\qquad
%\\
\widetilde{u}_y\big|_{y=-1}=-i\pi k \exp{(-|k|)},
%\end{array}\right 
$$
и, следовательно, задается формулой 
$$\widetilde{u}(k,y,\varepsilon)=-i\pi \sgn{(k)}\exp{(-|k|)}\left[\sgn{(k)}\ch\left(\sqrt{k^2-\varepsilon^2 k^4}(y+1)\right)+\right.$$
$$\left.+\frac{\sh\left(\sqrt{k^2-\varepsilon^2 k^4}(y+1)\right)}{\sqrt{1-\varepsilon^2 k^2}}\right].$$

Поэтому решение $u_{im}(x,y,\varepsilon)$ есть функция 
\begin{equation}\label{polfpl}
u_{im}(x,y,\varepsilon)=-i
\mathrm{Im}\int\limits_0^{\infty}
\frac{\exp{(\sqrt{k^2-\varepsilon^2k^4}(y+1))}}{2}\left(1+\frac{1}{\sqrt{1-\varepsilon^2k^2}}\right)
\exp{(-k-ikx)}dk-$$
$$-i\mathrm{Im}\int\limits_0^{\infty}
\frac{\exp{(-\sqrt{k^2-\varepsilon^2k^4}(y+1))
}}{2}\left(1-\frac{1}{\sqrt{1-\varepsilon^2k^2}}\right)
\exp{(-k-ikx)}dk.
\end{equation}
Раскладывая правые части формул (\ref{polfm}) и (\ref{polfpl}) в ряды Тейлора при $\varepsilon\to 0$ мы получаем представление эталонных решений $u_{re}(x,y,\varepsilon)$, $u_{im}(x,y,\varepsilon)$ в виде рядов теории возмущений (\ref{vnesh}), коэффициенты $U_{n,re}(x,y)$ и $U_{n,im}(x,y)$ которых по их главным членам -- функциям (\ref{rskri}) и (\ref{rskrr}) -- согласно вышесказанному задаются формулами 
%ИСПАВЛЕНИЕ 2
(\ref{rekk}):
%А НЕ
%(\ref{rekk}),(\ref{rskrn}):
 ввиду нарастающих вместе с номером $n$ особенностей его коэффициентов $U_{n}(x,y)$ при $x^2+y^2\to 0$ эти ряды теряет свою пригодность в малой окрестности начала координат. 

Для правильного описания асимптотики решения (\ref{polfm}) при малых $\varepsilon$ cделаем растяжение независимых переменных (\ref{skail}), которое дополним еще растяжением переменной интегрирования 
\begin{equation}\label{skaik}k=\frac{K}{\varepsilon^{2/3}}.\end{equation}  
После осуществления этих растяжений решение (\ref{polfm}) принимает вид  
$$\mathrm{Re}\int\limits_0^{\infty}
\frac{\exp{(\frac{K}{ \varepsilon^ {2/3}}
(\sqrt{1-\varepsilon^{2/3}K^2}(1+\varepsilon^{2/3}Y)-1)-iKX)}}
{2 \varepsilon^{2/3}}
\left(1+\frac{1}{\sqrt{1-\varepsilon^{2/3}K^2}}\right)dK
+$$
$$+%\underbrace{
\mathrm{Re}\int\limits_0^{\infty} \frac{\exp{(-\frac{K}{ \varepsilon^ {2/3}}
(\sqrt{1-\varepsilon^{2/3}K^2}(1+\varepsilon^{2/3}Y)+1)-iKX)}}
{2 \varepsilon^{2/3}}
\left(1-\frac{1}{\sqrt{1-\varepsilon^{2/3}K^2}}\right)dK.$$
Последний из интегралов, очевидно, есть ограниченная при $\varepsilon\to 0$ величина. 
Получаем, что правильное поведение эталонного решения (\ref{polfm}) в малой окрестности координат
задается соотношением 
\begin{equation}\label{polfmai} u(x,y,\varepsilon)=\frac{\mathrm{Re}I_0(X,Y,3)}{ \varepsilon^{2/3}}+\dots,\end{equation}
где $I_0(X,Y,3)$ есть интеграл Харди (\ref{pOrlv}) c индексами $m=3$ и $\nu=0$. 

Аналогично показывается, что соотношение 
\begin{equation}\label{polfmap} u(x,y,\varepsilon)=-i\frac{\mathrm{Im}I_0(X,Y,3)}{ \varepsilon^{2/3}}+\dots,\end{equation}
в главном по малому параметру $\varepsilon$ порядке описывает поведение при $x^2+y^2\to 0$ эталонного решения $u_{im}(x,y,\varepsilon)$, задаваемого формулой (\ref{polfpl}).

{\bf 2.3.} Рассмотрим теперь более общий случай чисто начальной задачи (\ref{sust}), (\ref{koshn}), про которую известно, что соответствующее решение $U_0(x,y)$ предельной задачи  Коши для уравнения Лапласа (\ref{lapl}), (\ref{Kof}) в начале координат имеет cингулярность, описываемую разложением (\ref{polm}). 

Возьмем разность $V(x,y)$ этой гармонической функции и точного решения 
%ИСПРАВЛЕНИЕ 3 
(\ref{rskri}).
%А НЕ 
%(\ref{rskrn}).
Гармоническая функция $V(x,y)$ уже {\it не имеет сигулярность} в начале координат. Поэтому в малой ее окрестности ряд теории возмущений (\ref{vnesh}), (\ref{rekk}) с функцией $V_0(x,y)$ в качестве главного члена сохраняет свою пригодность. А это значит, что в этой окрестности поведение решения $u(x,y,\varepsilon)$ задачи (\ref{sust}), (\ref{koshn}) с начальными данными $A(x)=U_0(x,-1)$ и $B(x)=(U_0(x,y))'_y|_{y=-1},$ определяемыми любым решением уравнения Лапласа $U_0(x,y)$, которое при $x^2+y^2\to 0$ разлагается в ряд (\ref{polm}), тоже в главном порядке задается соотношением (\ref{polfmai}).

Точно также показывается, что поведение в малой окрестности начала координат всех тех решений начальной задачи (\ref{sust}), (\ref{koshn}), для которых решения $U_0(x,y)$ предельной задачи Коши для уравнения Лапласа (\ref{lapl}), (\ref{Kof}) разлагаются при $x^2+y^2\to 0$ в ряд (\ref{polmm}), описывается соотношением (\ref{polfmap}). 

{\bf 2.4.} Модельный пример решения начально -краевой задачи (\ref{sust}), (\ref{koshn}), которое удовлетворяет условиям свободного опирания стержня (\ref{kra}) и описывается в окрестности начала координат посредством линейных комбинаций соотношений (\ref{polfmai}) и 
%ИСПРАВЛЕНИЕ 4
(\ref{polfmap})
%А НЕ
%(\ref{polfmai})
 дает сумма \begin{equation}\label{etb} u_{etb}(x,y,\varepsilon)=c_1u_{re}(x,y,\varepsilon)+c_2u_{im}(x,y,\varepsilon)+
u_{dir}(x,y,\varepsilon).\end{equation} 
Здесь $u_{re}(x,y,\varepsilon)$ и $u_{im}(x,y,\varepsilon)$ есть эталонные решения %ИСПРАВЛЕНИЕ 5
(\ref{polfm}) 
%А НЕ
%(\ref{polfpl})
и (\ref{polfpl}), а $u_{dir}(x,y,\varepsilon)$ решение уравнения (\ref{sust}), которое при $\varepsilon=0$ представляет собой решение $E(x,y)$ задачи Дирихле для уравнения (\ref{lapl}) в прямоугольнике $-2\leq y\leq q^2$, $-L\leq x\leq L$ c условиями по его краям 
$$E(x,y)|_{x=\pm L}= -c_1u_{re}(x,y,0)-c_2 u_{im}(x,y,0)$$ и согласованными с последними произвольными условиями на его верхней и нижней частях. 

Легко показать, что и произвольное $u(x,y,\varepsilon)$ решение начально-краевой задачи (\ref{sust}),(\ref{koshn}), (\ref{kra}), которое при $\varepsilon=0$ cводится к решению $U_0(x,y)$ 
%ИСПРАВЛЕНИЕ 6
уравнения Лапласа,
%А НЕ
%Лапласа,
 при $x^2+y^2\to 0$ представимого рядом 
$$U_0(x,y)=\frac{-c_1y}{y^2+x^2}+\frac{-c_2x}{y^2+x^2}+\dots,$$
также в окрестности начала координат описывается соотношением 
$$u(x,y,\varepsilon)=\frac{c_1\mathrm{Re} I_0(X,Y,3)}
{\varepsilon^{2/3}}-\frac{ic_2 \mathrm{Im}I_0(X,Y,3)}{\varepsilon^{2/3}}+\dots.$$
Для этого достаточно образовать разность 
$ W(x,y,\varepsilon)=u(x,y,\varepsilon)-u_{etb}(x,y,\varepsilon)$
этого решения с модельным решением (\ref{etb}) и заметить, что эта разность при $\varepsilon=0$ сводится к решению 
уравнения Лапласа, не имеющего особенности в начале координат.

{\bf 2.5.} Асимптотика интеграла $I_0(X,Y,3)$ при $Y\to-\infty$ легко вычисляется. После осуществления в нем замен 
\begin{equation}\label{minchp} X=Ys, \qquad K=\frac{\nu}{|Y|} \end{equation} 
получим, что для таких $Y$
\begin{equation}\label{minap}I_0(X,Y,3)=\frac{1}{|Y|}\int \limits_0^{\infty} \exp(-\nu(1-is)-\frac{\nu^3}{2|Y|^3})d\nu=$$
$$=\frac{1}{|Y|}\int\limits_0^\infty \exp(-\nu(1-is))\sum\limits_{m=0}^\infty\frac{1}{m!}\left(-\frac{\nu^3}{2|Y|^3}\right)^md\nu=$$%\frac1\l\sum\limits_{m=0}^\infty\frac{(3m)!}{m!}\left(\frac{-1}{2|Y|^3}\right)^3(1-is)^{-3m-1}=$$
$$=\frac{1}{-Y}\sum\limits_{m=0}^\infty\frac{(3m)!}{m!}\left(\frac{1}{2Y^3}\right)^m(1-is)^{-3m-1}=\frac{1}{-Y+iX}+\dots.\end{equation}
С помощью рассуждений, подобных тем, что были использованы при описании асимптотики (\ref{minap}), нетрудно показать, что на самом деле 
и всюду вне сектора 
\begin{equation}\label{Sek} Y>0, \qquad |X|<\sqrt{3}Y.\end{equation}
 поведение интеграла $I_0(X,Y,3)$ при $X^2+Y^2\to \infty $ также задается правой частью соотношения (\ref{minap}). Справедливость этой асимптотики означает, что асимптотики (\ref{polfmai}), (\ref{polfmap}) в главном порядке согласуются (сшиваются) с рядами теории возмущений (\ref{vnesh}), (\ref{rekk}), задаваемых функциями (\ref{polm}) и, соответственно, (\ref{polmm}). 

В полуплоскости $Y \geq 0$ полное и равномерно пригодное асимптотическое разложение этого интеграла при $X^2+Y^2\to \infty $ было описано в работе \cite{Back}. В главном порядке данное разложение задается формулой ($|\arg(Y-iX)|\leq \pi/2$)
$$I_0(X,Y,3)=\frac{1}{-Y+iX}\left(1+O\left(\frac{1}{|X+iY|^3}\right)\right)+$$
$$+\sqrt{\pi}\left(\frac{2}{3(Y-iX)}\right)^{1/4}\exp {\left(\left(\frac{2}{3}(Y-iX)\right)^{3/2}\right)}
\left(1+O(|X+iY|^{-3/2})\right),$$
% \begin{equation}\label{Back} %I_0(X,Y,3)=\frac{1}{-Y+iX}\left(1+O\left(\frac{1}{|X+iY|^3}\right)\right)+$$
%$$+\sqrt{\pi}\left(\frac{2}{3(Y-iX)}\right)^{1/4}\exp %{\left(\left(\frac{2}{3}(Y-iX)\right)^{3/2}\right)}
%\left(1+O(|X+iY|^{-3/2})\right), \end{equation} которое, в частности, означает, что внутри сектора (\ref{Sek}) модуль интеграла $I_0(X,Y,3)$ при $Y\to \infty$ экспоненциально растет. 

{\bf 2.6.} Перейдем к конкретному примеру такой задачи. В нем решение предельной задачи задается рядом Фурье 
$$U_0(x,y)=-\exp{(y)}\sin{x}-\exp{(2y)}\sin{2x}-\sum \limits_{j=3}^{\infty} \exp{(jy)}\sin{jx}.$$ 
Это есть \cite[формула 5.4.9(3)]{Tolst} ряд Фурье функции
\begin{equation}\label{tolst}U_0(x,y)=-\frac{\exp{(y)}\sin{x}}{1-2 \exp{(y)} 
\cos{x}+\exp{(2y)}},\end{equation}
которая бесконечно дифференцируема при всех отрицательных $y$. В этой полуплоскости данная функция, очевидно, удовлетворяет 
%ИСПРАВЛЕНИЕ 7 
уравнению 
%А НЕ 
%уравнения 
Лапласа (\ref{lapl}) и краевым условиям $u_0(x,y)|_{x=\pm \pi}=u_0(x,y)''_{xx}|_{x=\pm \pi}=0$. В начале координат она имеет особенность, описываемую асимптотикой (\ref{polmm}).

Ряд теории возмущений (\ref{vnesh}), (\ref{rekk}), 
определяемый этой функций $U_0(x,y)$, за исключением малой окрестности начала координат описывает в прямоугольнике $-1\leq y\leq 0$, $-\pi\leq x\leq \pi$ поведение решения $u(x,y,\varepsilon)$ уравнения (\ref{sust}), которое удовлетворяет начальным условиям 
%\begin{equation}\label{tolstn}u(x,-1)= -\frac{\exp{(1)}\sin{x}}{1-2 \exp{(1)} 
%\cos{x}+\exp{(2)}}, u'_y(x,-1)=\frac{(\exp(1)-\exp{(3)})\sin{x}}{(1-2 \exp{(1)} 
%\cos{x}+\exp{(2)})^2}\end{equation}
\begin{equation}\label{tolstn}u(x,-1)= -\frac{e\sin{x}}{1-2 e 
\cos{x}+e^2},\qquad u'_y(x,-1)=\frac{(e-e^3)\sin{x}}{(1-2 e 
\cos{x}+e^2)^2}\end{equation}
и краевым условиям свободного опирания 
\begin{equation}\label{tolstb}u(x,y)|_{x=\pm \pi}=u(x,y)''_{xx}|_{x=\pm \pi}=0.\end{equation}
Но понятно, что при малых $x$ и $y$ главный член этого ряда $U_0(x,y)$ имеет особенность, описываемую дифференцируемой асимптотикой (\ref{polmm}). Поэтому истинное поведение решения $u(x,y,\varepsilon)$ уравнения (\ref{sust}), удовлетворяющего условиям (\ref{tolstn}) и (\ref{tolstb}), в малой окрестности начала координат описывается не этим рядом, а соотношением (\ref{polfmap}). Данные факты хорошо иллюстрируют рисунки 1 и 2. На рисунке 3 также четко видно, что именно в малой окрестности начала координат, являющейся точкой сингулярности решения $U_0(x,y)$ предельной задачи, действительно расположена расширяющаяся по мере роста $y$ область быстро нарастающих значений амплитуды $u(x,y,\varepsilon)$.

\begin{figure}[htp]
    \centering
        \includegraphics[width=80mm,height=60mm]{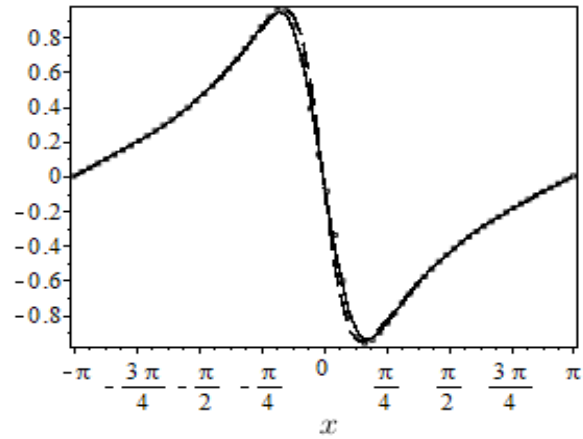}
    \caption{Функции $u(x,-0.5,\varepsilon)$ (сплошная линия), $u_0(x,-0.5,\varepsilon)$ (штриховая линия) и $u_0(x,-0.5,\varepsilon)+\varepsilon^2u_1(x,-0.5,\varepsilon)$ 
		(пунктирная линия) при $\varepsilon=0.1$.}
    \label{fig6}
\end{figure}

\begin{figure}[htp]
    \centering
        \includegraphics[width=80mm,height=60mm]{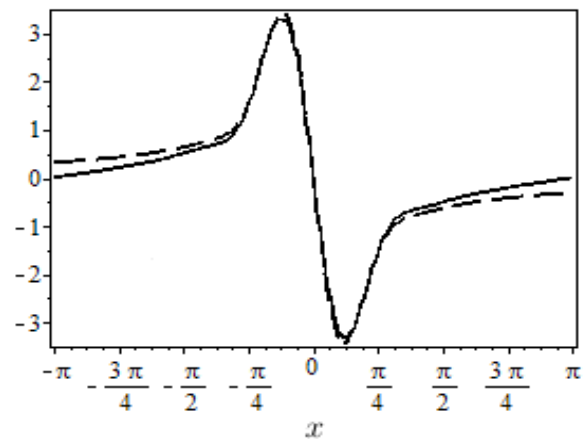}
    \caption{Функции $u(x,0,\varepsilon)$ (сплошная линия) и $\varepsilon^{-2/3}\mathrm{Im} I_0(x\varepsilon^{-2/3},0)$ при $\varepsilon=0.1$.}
    \label{fig7}
\end{figure}

\begin{figure}[htp]
    \centering
        \includegraphics[width=80mm,height=70mm]{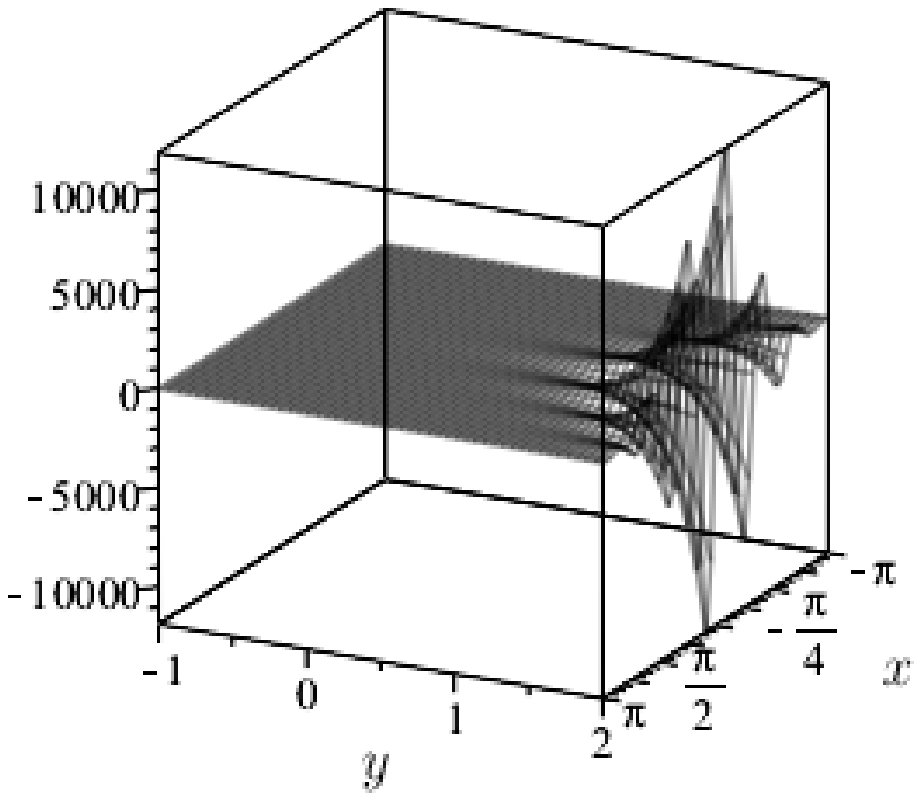}
    \caption{Функциия $u(x,y,\varepsilon)$  при $\varepsilon=0.1$.}
    \label{fig8}
\end{figure}

\section{КОРНЕВАЯ СИНГУЛЯРНОСТЬ }

{\bf 3.1.} Типичные особенности решений уравнений Лапласа (\ref{lapl})
, которые в главном порядке задаются вещественной частью корня (\ref{root}), при $x^2+y^2\to 0$ имеет 
дифференцируемые по $x$ и $y$ асимптотики 
\begin{equation} \label{rskr}U_0(x,y)
=\sqrt{-y+\sqrt{y^2+x^2}}+\dots.\end{equation} 
%Характер перестроек такого сорта решений, которое происходит сними в окретсности начала координат, во многом напоминает перестройки, описанные в предыдущем разделе. Поэтому мы часто будем ссылаться на уже сделанные там описания. Правда, перестройки , рассматривамые ниже, потребовали все же несколько более сложной техники.  

По аналогии с предыдущим разделом логично прежде всего рассмотреть эталонное решение $u_{sq}(x,y,\varepsilon)$ чисто начальной задачи (\ref{sust}), (\ref{koshn}) с условиями 
$$u_{sq}|_{y=-1}=\sqrt{1+\sqrt{1+x^2}}, \qquad (u_{sq})'_y|_{y=-1}=\frac{1-\sqrt{1+x^2}}{2\sqrt{1+\sqrt{1+x^2}}\sqrt{1+x^2}},$$
%\begin{equation}\label{koshs}u_{sq}|_{y=-1}=\sqrt{1+\sqrt{1+x^2}}, \qquad %(u_{sq})'_y|_{y=-1}=\frac{1-\sqrt{1+x^2}}{2\sqrt{1+\sqrt{1+x^2}}\sqrt{1+x^2}},\end{equation}
определяемыми частным случаем 
\begin{equation} \label{rskrn}U_0(x,y)=Q(x,y)
\equiv\sqrt{-y+\sqrt{y^2+x^2}}\end{equation} 
решений уравнения Лапласа с асимптотиками (\ref{rskr}).

{\bf 3.2.} Технически удобнее, однако, предварительно описать его вторую производную $u_{sq}(x,y,\varepsilon)''_{xx}$, являющуюся решением начальной задачи (\ref{sust}), (\ref{koshn}), 
где 
\begin{equation} \label{koshABn} A(x)=\frac{\partial^2}{\partial x^2}\sqrt{1+\sqrt{1+x^2}}, \qquad B(x)=-\frac{\partial^2}{\partial x^2}\frac{\sqrt{1+\sqrt{1+x^2}}}{2\sqrt{1+x^2}}.\end{equation}

Из четности начальных данных (\ref{koshABn}) по переменной $x$ cледует четность их преобразований Фурье (\ref{inpf}) $\widetilde{A}(k)$ и $\widetilde{B}(k)$ по переменной $k$. Поэтому из формулы $$\frac{\cos(\frac32\arctg x)}{2\sqrt{2}(x^2+1)^{3/4}}=(Q(x,-1))''_{xx}=A(x)$$ и известного \cite[формула 1.4(3)]{Beit} косинус -- преобразования Фурье 
$$\int \limits_{0}^\infty k^{1/2}\exp{(-k)}\cos(kx)dk=\frac{\sqrt{\pi}\cos(\frac32\arctg x)}{2(x^2+1)^{3/4}} $$ cледует, что 
%\begin{equation}\label{ankf}
$$\widetilde{A}(k)= \sqrt{\frac{\pi}{2}}|k|^{1/2}\exp{(-|k|)}.$$
%\end{equation}
И поскольку $$ (Q_(x,y))'_y|_{y=-1}=-\frac12\sqrt{1+\sqrt{x^2+1}}
+x\left(\sqrt{1+\sqrt{x^2+1}}\,\right)'_x,$$ то
$$\widetilde{B}(k)=\int\limits_{-\infty}^\infty (Q(x,y))'''_{xxy}|_{y=-1}\exp{(ikx)} dx
%Решая эту начальную задачу и, применяя затем обратное преобразование Фурье, находим, что 
=\sqrt{\frac{\pi}{2}}|k|^{3/2}\exp{(-|k|)}.$$
Таким образом, после преобразования Фурье задача (\ref{sust}), (\ref{koshn}), (\ref{koshABn}) сводится к решению задаче Коши для уравнения (\ref{oduf}) с начальными 
с условиями 
$$\widetilde{u}\big|_{y=-1}=\sqrt{\frac{\pi}{2}}|k|^{1/2}\exp{(-|k|)},\qquad
%\\
\widetilde{u}_y\big|_{y=-1}=\sqrt{\frac{\pi}{2}}|k|^{3/2}\exp{(-|k|)}.$$
% =\sqrt{\frac{\pi}{2}}|k|^{3/2}\exp{(-|k|)}.$$

Решая ее и, применяя затем обратное преобразование Фурье, находим, что %=\sqrt{\frac{\pi}{2}}|k|^{3/2}\exp{(-|k|)}.$$
\begin{equation}\label{susk} u_{sq}(x,y,\varepsilon)''_{xx}=\frac{1}{2\pi}\int\limits_{-\infty}^\infty
\widetilde{u}(k,y)\exp{(-ikx)}dk=\frac{1}{\pi}\mathrm{Re}\int\limits_{0}^\infty
\widetilde{u}(k,y)\exp{(-ikx)}dk=$$
%$$=\frac{1}{\sqrt{2\pi}}\mathrm{Re}
%\int\limits_0^\infty k^{1/2}\exp{(-k)}\exp{(-ikx)}
%\Big(\ch(\sqrt{k^2-\varepsilon^2k^4}(y+1))+$$ %$$+\frac{\sh(\sqrt{k^2-\varepsilon^2k^4}(y+1))}{\sqrt{1-\varepsilon^2k^2}}\Big)dk=$$
$$=\frac{1}{2\sqrt{2\pi}}\mathrm{Re}\bigg[
\int\limits_0^\infty k^{1/2}\exp{(-k-ikx)}
\exp{(\sqrt{k^2-\varepsilon^2k^4}(y+1))}\Big(1+\frac{1}{\sqrt{1-\varepsilon^2k^2}}\Big)dk+$$
$$+
\int\limits_0^\infty k^{1/2}\exp{(-k-ikx)} \exp{(-\sqrt{k^2-\varepsilon^2k^4}(y+1))}\Big(1-\frac{1}{\sqrt{1-\varepsilon^2k^2}}\Big)
\bigg]dk. \end{equation}

Почти слово в слово повторяя теперь рассуждения п. 2.2, которые приведены там после формулы 
(\ref{polfpl}), приходим к выводу:

%-- ряд теории возмущений (\ref{vnesh}), (\ref{rekk}), 
%в котором главный член есть вторая производная по $x$ точного решения (\ref{rskrn}) уравнения Лапласа правильно задает поведение $u_{sq}(x,y,\varepsilon)''_{xx}$ при $\varepsilon \to 0$
%при $-1\leq y\leq 0$ за исключением малой окрестности начала координат; 

-- в малой окрестности начала координат правильное поведение $u_{sq}(x,y,\varepsilon)''_{xx}$ задается соотношением 
\begin{equation}\label{int}u_{sq}(x,y,\varepsilon)''_{xx}=\frac{\mathrm{Re}I_{1/2}(X,Y,3) }{\varepsilon \sqrt{2\pi}}+\dots,\end{equation}
где $I_{1/2}(X,Y,3)$ есть интеграл Харди (\ref{pOrlv}) c $m=3$, $\nu=1/2$ и переменными 
$X, Y$, заданными растяжениями (\ref{skail}).

{\bf 3.3.} Из сказанного выше ясно, что ряд теории возмущений (\ref{vnesh}), (\ref{rekk}), главный член 
%ИСПРАВЛЕНИЕ 8
$U_0(x,y)$ 
%А НЕ 
%$U_0(x,y
которого есть точное 
%ИСПРАВЛЕНИЕ 9
решение 
%А НЕ
%решения 
(\ref{rskrn}) уравнения Лапласа (\ref{lapl}), правильно задает поведение эталонного решения $u_{sq}(x,y,\varepsilon)$ уравнения (\ref{sust}) при $\varepsilon \to 0$
при $-1\leq y\leq 0$ за исключением малой окрестности начала координат. 

В свою очередь, из этого факта и из соотношения (\ref{int}) cледует, что в этой малой окрестности координат поведение эталонного решения $u_{sq}(x,y,\varepsilon)$ в главном по малому параметру $\varepsilon$ порядке задается соотношением 
\begin{equation}\label{in}u(x,y,\varepsilon)=-\frac{(\varepsilon)^{1/3}\mathrm{Re}[2(Y-iX) I_{-1/2}(X,Y,3)-3I_{3/2}(X,Y,3)]}{\sqrt{2\pi}}+\dots,\end{equation}
где $X$ и $Y$ есть растянутые переменные (\ref{skail}), а $I_{-1/2}(X,Y,3)$ и $I_{3/2}(X,Y,3)$
представляют собой частные случаи интеграла Харди (\ref{pOrlv}).

Действительно, из вида главного члена (\ref{rskrn}) ряда теории возмущений 
(\ref{vnesh}), (\ref{rekk}) и из соображений, изложенных в последних двух абзацах раздела 1 настоящей статьи, вытекает, что это поведение описывается соотношением 
$$u_{sq}(x,y,\varepsilon)=\varepsilon^{1/3}u_{in}(X,Y)+\dots. $$
где $u_{in}(X,Y)$ есть точное решение уравнения Лапласа (\ref{Lapl}), которое при $-Y+X^2\to\infty$ имеет алгебраическую асимптотику 
\begin{equation} \label{alg}u_{in}(X,Y)=\sqrt{-Y+\sqrt{Y^2+X^2}}+\dots .\end{equation}
Из cоотношения (\ref{int}) выводим, что эта гармоническая функция (очевидно, четная по переменной $X$) необходимо имеет вид 
$$u_{in}(X,Y)=-\frac{\mathrm{Re}\int\limits_{0}^{\infty}k^{-3/2}
[\exp{(k(Y-iX)-k^3/2)}-1] dk}{\sqrt{2\pi}}+C_1Y+C_2=$$
$$=-\frac{\mathrm{Re}[2(Y-iX) I_{-1/2}(X,Y,3)-3I_{3/2}(X,Y,3)]}{\sqrt{2\pi}}+C_1Y+C_2$$
%\begin{equation}\label{promuin}u_{in}(X,Y)=-\frac{\mathrm{Re}\int\limits_{0}^{\infty}%k^{-3/2}
%[\exp{(k(Y-iX)-k^3/2)}-1] dk}{\sqrt{2\pi}}+C_1Y+C_2=$$
%$$=-\frac{\mathrm{Re}[2(Y-iX) %I_{-1/2}(X,Y,3)-3I_{3/2}(X,Y,3)]}{\sqrt{2\pi}}+C_1Y+C_2\end{equation}
где $C_1$ и $C_2$ -- вещественные постоянные.
Справедливость соотношения (\ref{in}), а также соотношения (\ref{alg}) (в частности, равенство нулю постоянных $C_1$ и $C_2$) следует теперь 
из известных \cite[формула (21.18)]{Rie} асимптотик интеграла Харди (\ref{pOrlv}) 
 при $-Y+X^2\to\infty$.

 Более того, из этих асимптотик интеграла Харди следует, что при $Y^2+X^2\to\infty$ соотношение (\ref{alg}) cправедливо и всюду вне сектора (\ref{Sek}). А вот внутри этого сектора функция $u_{in}(X,Y)$ при больших значениях аргументов уже экспоненциально растет --  в главном порядке ее асимптотика при $Y^2+X^2\to \infty$ внутри данного сектора описывается формулой: 
%$$u_{in}(X,Y)=\mathrm{Re}\frac{ %\sqrt{3}}{2(Y-iX)}\exp{[(\frac{2(Y-iX)}{3})^{3/2}]}\left(1+O(|Y-iX|^{-3/2})\right)$$
\begin{equation}\label{exa} u_{in}(X,Y)=\mathrm{Re}\frac{ \sqrt{3}}{2(Y-iX)}\exp{
\Big[\Big(\frac{2(Y-iX)}{3}\Big)^{3/2}\Big]}\left(1+O(|Y-iX|^{-3/2})\right)\end{equation}
(с помощью интегрирования по частям  эта асимптотика выводится из равенства 
$$u_{in}(X,Y)'_Y=\mathrm{Re}\frac{I_{-1/2}(X,Y,3)}{\sqrt{2\pi}},$$ получающейся в результате дифференцирования по переменной $Y$ формулы (\ref{in}),
и из асимптотик в этом секторе интегралов $I_{-1/2}(X,Y,3)$, $I_{3/2}(X,Y,3)$). 

{\it Замечание.} В экспоненциальной части асимптотики интеграла Харди, описываемой формулой 
(21.18) монографии \cite{Rie}, имеется неточность -- верный ответ (\ref{exa}) выписан после устранения этой неточности. (Применяя так называемую лемму Ватсона, автор \cite{Rie} ошибочно опустил деление на множитель $\beta$, который в нашем случае равен 2.)

{\bf3.4.} Поведение любого решения чисто начальной задачи $u(x,y,\varepsilon)$ задачи (\ref{sust}), (\ref{koshn}) с начальными данными $A(x)=U_0(x,-1)$ и $B(x)=(U_0(x,y))'_y|_{y=-1},$ задаваемыми решением $U_0(x,y)$ уравнения Лапласа (\ref{lapl})  с асимптотикой (\ref{rskr}) при $x^2+y^2\to 0$, в малой окрестности начала координат также описывается соотношением (\ref{in}). 

Чтобы убедиться в справедливости последнего утверждения, аналогично п. {\bf2.3} рассмотрим разность $R(x,y)=U_0(x,y)-\sqrt{-y+\sqrt{y^2+x^2}}$ этого решения уравнения Лапласа и правой части формулы (\ref{rskrn}). В отличие от ситуации, описанной в разделе п. {\bf2.3}, функция $R(x,y)$ не будет, вообще говоря, бесконечно дифференцируема в начале координат. Однако справедливость асимптотики (\ref{rskr}) означает, что порядок ее сингулярности в начале координат {\it меньше}, чем у составляющих этой разности. Поэтому, хотя ряд теории возмущений (\ref{vnesh}),(\ref{rekk}) с гармонической функцией $R(x,y)$ в качестве главного члена в ситуации общего положения и теряет свою пригодность при $x^2+y^2\to0$, но все же соотношение (\ref{in}), действительно, по-прежнему справедливо.

{\bf 3.5.} Переходя к решениям начально-краевых задач, удовлетворяющих условиям свободного опирания (\ref{kra}), прежде всего сразу оговорим то, что предельными для рассматриваемых здесь будут решения $U_{0}(x,y)$ уравнения Лапласа (\ref{lapl}), которые, вообще говоря, в начале координат в нуль не обращаются. Предполагается, что в ее малой окрестности дифференцируемые асимптотики, задаваемые правой частью соотношения (\ref{rskr}), справедливы лишь для отклонений $U_{0}(x,y)$ от их значений в нуле $U_0(0,0)$, то есть, что при $x^2+y^2\to 0$ 
$$U_0(X,Y)=U_0(0,0)+\sqrt{-y+\sqrt{y^2+x^2}}+\dots.$$
%\begin{equation}\label{difero} %U_0(X,Y)=U_0(0,0)+\sqrt{-y+\sqrt{y^2+x^2}}+\dots.\end{equation}

Приведем пример модельный пример $u_{msk}(x,y,\varepsilon;L)$ такого решения начально -краевой задачи для уравнения (\ref{sust}). Для этого, действуя аналогично тому, как в п. {\bf 2.3}, введем в рассмотрение функцию $H(x,y;L)$, являющуюся решением задачи Дирихле
$$ H''_{xx}+H''_{yy}=0,\quad -L\leq x\leq L,\;-3\leq y<1,$$
$$ H\Big|_{x=\mp L}=Q(\mp L,y)=\sqrt{-y+\sqrt{y^2+L^2}},$$
$$ H\Big|_{y=-3}=Q(x,-3)=\sqrt{3+\sqrt{9+x^2}},$$
$$ H\Big|_{y=1}=
\sqrt{3+\sqrt{9+x^2}}+\sqrt{-1+\sqrt{1+L^2}}-\sqrt{3+\sqrt{9+L^2}}$$
%$$\left\{
%\begin{array}{l}
%\ds \Delta H=0,\quad -L<x<L,\;-3<y<1,\\
%\ds H\Big|_{x=\mp L}=Q(\mp L,y)=\sqrt{-y+\sqrt{y^2+L^2}},\\
%\ds H\Big|_{y=-3}=Q(x,-3)=\sqrt{3+\sqrt{9+x^2}},\\
%\ds H\Big|_{y=1}=Q(x,-3)+Q(L,1)-Q(L,-3)=\\
%\sqrt{3+\sqrt{9+x^2}}+\sqrt{-1+\sqrt{1+L^2}}-\sqrt{3+\sqrt{9+L^2}}
%\end{array}
%\right.$$

Образуем затем разность $S(x,y;L)=\sqrt{-y+\sqrt{y^2+x^2}}-H(x,y;L)$ и построим по ней модельное решение $u_{msk}(x,y,\varepsilon;L)$ уравнения (\ref{sust}), которое удовлетворяет условиям свободного опирания (\ref{kra}) и начальным условиям 
%\begin{equation}\label{mosk}
$$u_{msk}|_{y=-1}=S(x,-1;L),\qquad (u_{msk})'_y|_{y=-1}=S(x,y;L)'_y|_{y=-1}.$$
%\end{equatoion}

Графики функций $Q(x,y)$, %$H(x,y,3)$, 
$S(x,y;3)$ и $u_{msk}(x,y, \varepsilon;3)$ изображены на рисунке 4 и 5. Из последнего из этих графиков мы видим, что в окрестности начала координат начинается расширяющаяся область нарастающих значений $|u_{msk}(x,y, \varepsilon;3)|$. Рисунок 6 иллюстрирует то обстоятельство, что ряд теории возмущений (\ref{vnesh}), (\ref{rekk}) с главным членом $U_0(x,y)=S(x,y;3)$ правильно приближает модельное решение $u_{msk}(x,y, \varepsilon)$ при $y=-0.7$. 

\begin{figure}[htp]
    \centering
        \includegraphics[width=80mm,height=65mm]{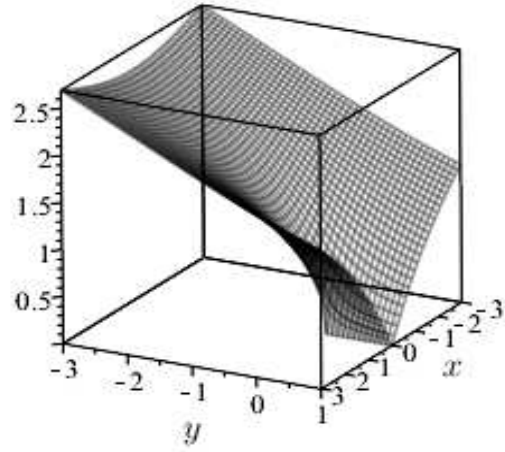}
    \caption{Функции $Q(x,y)=\sqrt{-y+\sqrt{x^2+y^2}}$.}
    \label{fig8}
\end{figure}

%\begin{figure}[htp]
%    \centering
%        \includegraphics[width=80mm,height=65mm]{H.eps}
%    \caption{График функции $H(x.y;3)$.}
%    \label{fig9}
%\end{figure}

\begin{figure}[htp]
    \centering
        \includegraphics[width=80mm,height=65mm]{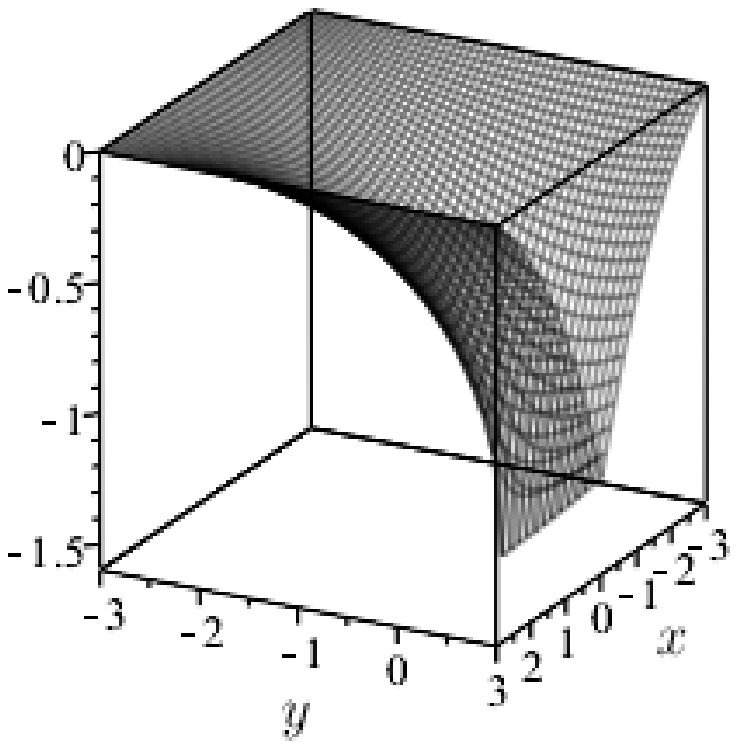}
    \caption{Функции $S(x.y;3)$.}
    \label{fig10}
\end{figure}

\begin{figure}[htp]
    \centering
        \includegraphics[width=80mm,height=65mm]{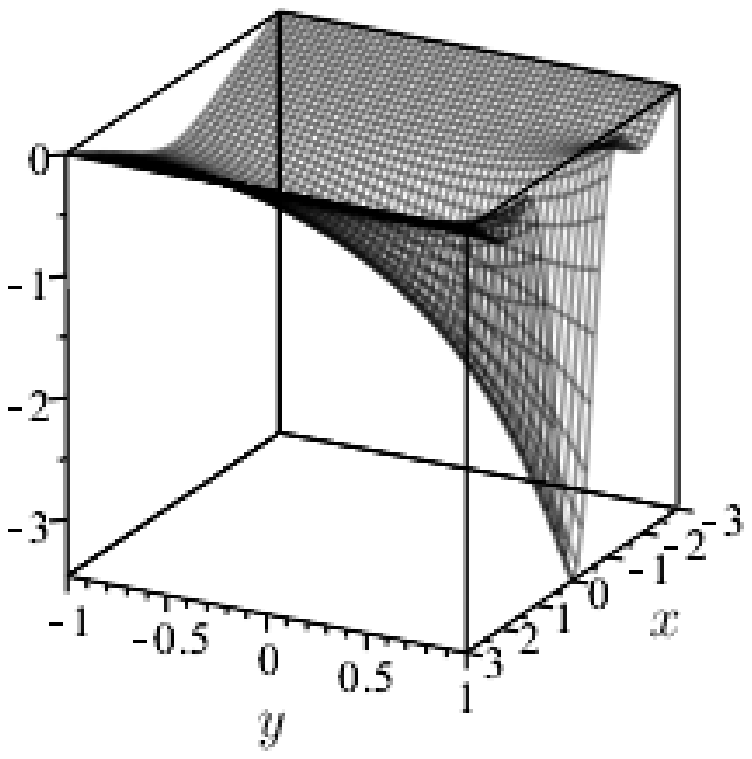}
    \caption{Функция $u_{msk}(x,y,\varepsilon;3)$  при $\varepsilon=0.1$.
}
    \label{fig11}
\end{figure}

\begin{figure}[htp]
    \centering
        \includegraphics[width=80mm,height=60mm]{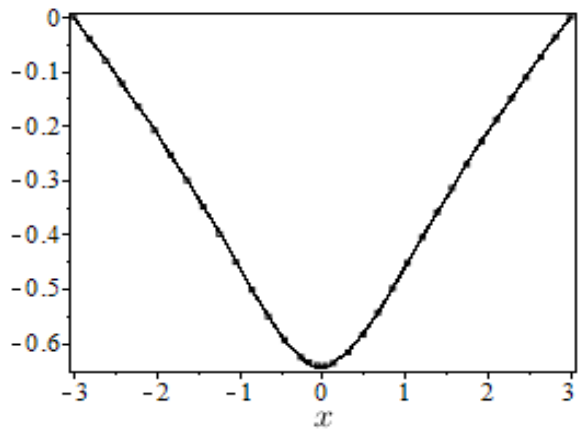}
    \caption{Функции $u_{msk}(x,-0.7,0.1;3)$ (сплошная линия) и $S(x,y;3)$.}
    \label{fig12}
\end{figure}

Из рассмотрения этого модельного решения при помощи рассуждений рассуждений, опять-таки почти полностью повторяющих приведенные в разделе {\bf 2.3}, приходим к выводу: 
решения начально-краевых задач (\ref{sust}),(\ref{kra}) с начальными данными
$$u|_{y=-1}=U_0(x,-1), \qquad u'_y|_{y=-1}=(U_0(x,y))'_y|_{y=-1},$$
которые определяются решениями $U_0(x,y)$ предельного уравнения Лапласа, описанными в первом абзаце данного пункта, в малой окрестности начала координат в главном по малому параметру $\varepsilon$ порядке описываются соотношением 
$$u(x,y,\varepsilon)=U_{0}(0,0)-\frac{(\varepsilon)^{1/3}\mathrm{Re}[2(Y-iX) I_{-1/2}(X,Y,3)-3I_{3/2}(X,Y,3)]}{\sqrt{2\pi}}+\dots. $$

Рисунок 8 иллюстрирует справедливость этого соотношения для модельного решения $u_{msk}(x,y,\varepsilon)$ при $Y=0$.

\begin{figure}[htp]
    \centering
        \includegraphics[width=80mm,height=60mm]{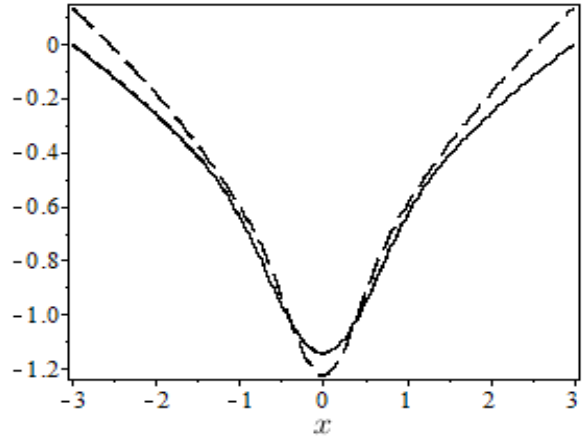}
    \caption{Функции $u(x,y,\varepsilon)$ (сплошная линия) и 
$-H_(0,0)+\varepsilon^{1/3}\mathrm{Re}[2(Y-iX) I_{-1/2}(X,Y;3)+3 I_{-1/2}(X,Y;3)] (2 \pi)^{-1/2}$
при $y=0$ и $\varepsilon=0.1$.}
    \label{fig13}
\end{figure}

%Поскольку $H(x,y;L)$ --- гладкая и ограниченная, то фу%нкция $S(x,y;L)$ 
%имеет ту же особенность в нуле, что и $U_0(x,y)$, но в то же время является
%главным членом асимптотики решения задачи (при $x<0$)\vspace{-2mm}:
%$$\left\{
%\begin{array}{l}
%v_{xx}+v_{yy}+\ep^2 v_{xxxx}=0,\quad -\infty<x<\infty,\;y>-1,\\
%v\big|_{x=-L}=v\big|_{x=a}=v_{xx}\big|_{x=-L}=v_{xx}\big|_{x=a}=0,\\
%v\big|_{y=-1}=S(x,-1),\\
%v_y\big|_{y=-1}=S_y(x,-1)
%=U_{0y}(x,-1)-H_y(x,-1;L)=0.
%\end{array}\right.$$

\section{Заключение}

Изучению вопроса о сглаживающем влиянии членов дифференциальных уравнений с малым параметром при старших производных на типичные сингулярности решений предельных {\it нелинейных} уравнений посвящено большое количество самых различных исследований. %(влияние эффектов малой вязкости или дисперсии на особенности решений гидродинамических уравнений и др.). 
Это влияние описывалось и для ситуаций, когда предельные задачи линейны, однако соответствующие им линейные дифференциальные уравнения имеют сингулярности в коэффициентах. Хорошо известны, в частности, результаты об истинном поведении решений линейных дифференциальных уравнений с малым параметром при производных в окрестностях каустик решений системы уравнений приближения ВКБ, одни из которых -- уравнения эйконала -- нелинейны, а коэффициенты других -- линейных уравнений переноса -- определяются решениями уравнений эйконала и имеют особенности в каустических точках.

Как ни странно, но по-видимому до настоящей работы этот естественный вопрос для чисто линейных задач, предельными к которым являются также {\it линейные} задачи с сингулярностями их решений {\it внутри} рассматриваемой области, не связанными с особенностями коэффициентов соответствующего предельного линейного дифференциального уравнения, ни в одной из публикаций не подымался. Между тем, очевидна обширность множества таких задач. Описание особенностей изгибов стержня при его сильном продольном сжатии, представленное в данной работе, является, вероятно, первым из исследований этого класса проблем. 

За оказавшуюся весьма полезной при написании данной статьи информацию мы выражаем свою признательность Н. Ф. Валееву, А. Р. Данилину, М. А. Ильгамову и И. В. Мельниковой. Работа первого из авторов над статьей выполнена за счет гранта РНФ
(проект 15-11-10018). %и 14-11-00078)
%\newpage

\def\refname{\centerline{\underline{\hspace{3 cm }}}}


\begin{thebibliography}{99}


\bibitem{Vol}% А. С. Вольмир, {\it Устойчивость деформируемых систем} (Наука,1967). 
 A. S. Vol’mir, {\it Stability of Deformable Systems} (Nauka, 1967) [in russian].

\bibitem{PaG}% Я. Г. Пановко, И. И. Губанова, {\it Устойчивость и колебания упругих систем}, (Наука, 1987). 
 Ya. G. Panovko and I. I. Gubanov, {\it Stability and Vibrations of Elastic Systems} (Fizmatlit, 1987) [in
russian].



\bibitem{LaI} %М. А. Лаврентьев, А. Ю. Ишлинский, %Дина-мические формы потери устойчивости упругих си-стем  
%ДАН СССР. {\bf64}, 779 (1949).
 M. A. Lavrent’ev and A. Yu. Ishlinskii, ``Dynamic buckling modes of elastic systems,'' Dokl. Akad.
Nauk 64 (6), 779--782 (1949).



\bibitem{Il} %М. А. Ильгамов, {\it Колебания упругих оболочек, содержащих жидкость и газ}, (Наука, 1969). 
M. A. Il’gamov, {\it Vibrations of Elastic Shells Containing a Liquid and Gas} [in Russian], (Nauka, 1969) [in russian].

\bibitem{Pa} M. P. Paidouisis, ``Flatter of conservative sysytems of pipes conveying incompressible fluid,''
J. Mech. Ing. Sci.  {\bf17}, 19--25 (1975).

\bibitem{Vo} %А. С. Вольмир, {\it Оболочки в потоке жидкости и газа} (Наука, 1979). 
 A. S. Vol'mir{\it Shells in Fluid Flow. Hydroelastic Problems}, 
 
(Nauka, 1979).
\bibitem{Ah} %А. М. Ахтямов, Г. Ф. Сафина, %Определение виброзащитности трубопровода
%Пр. мех. и тех. физ. {\bf49}, 139 (2008).
A. M. Akhtyamov, G. F. Safina, ``Vibration-proof conduit,'' Journal of Applied Mechanics and Technical Physics, {\bf49} (1), 114--121 (2008).


\bibitem{s46}% А. Б. Шварцбург, {\it Геометрическая оптика в нелинейной теории волн} (Наука, 1976).
 A. B. Shvartsburg, {\it Geometrical optics in the Nonlinear Theory of Waves} (Nauka, 1976) [in russian].

\bibitem{ZhT}% С. Л. Жданов, Б. А. Трубников, {\it Квазигазовые неустойчивые среды} (Наука, 1991) 
 S. K. Zhdanov, B. A. Trubnikov, {\it Quasigaseous Unstable Media}
(Nauka, 1991) [in russian].


\bibitem{Gel'fSh} %И. М. Гельфанд, Г. Е. Шилов, {\it Обобщенные функции. Выпуск 3. Некоторые вопросы теории дифференциальных уравнений} (Государственное издательство физико-математической литературы, 1958).
I. M. Gel’fand, G. E. Shilov, {\it Generalized Functions, Vol. 3, Some Questions of the Theory of Differential Equations} (Academic Press, 1968).


\bibitem{Isch} А. Ю. Ишлинский, {\it Прикладные задачи механики. Книга 2. Механика упругих и абсолютно твердых тел} (Наука, 1976).

\bibitem{Ishl} А. Ю. Ишлинский, {\it  Механика, идеи, задачи, приложения} (Наука, 1985).

\bibitem{Kon} C. Koning, J. Taub, ``Impact buckling of thin bars in the elastic range hinged
at both ends"', Luftfahrtforshung {\bf 10} (2), 55--64 (1933); 
also transl. as NACA TM 748 (1934).

\bibitem{Mai}% М. А. Ильгамов, ДАН,  {\bf 432}, 624 (2010).
M. A. Il'gamov, ``Reconstruction of harmonics during the dynamic loss of stability in mechanical systems,'' Doklady Physics, {\bf55} (6), 297–301 (2010).

\bibitem{Mail}% М. А. Ильгамов, ДАН,  {\bf 457}, 656 (2014).
M. A. Il'gamov, ``Dependence of dynamic buckling of a rod on the initial conditions,'' 
Doklady Physics, {\bf59} (8), 385--388 (2014).

\bibitem{Mor}% Н. Ф. Морозов, А. К. Беляев, П. Е. Товстик, Т. П. Товстик,  ДАН,  {\bf 463}, 543 (2014).
N. F. Morozov, A. K. Belyaev, P. E. Tovstik, T. P. Tovstik, ``The Ishlinskii—Lavrent’ev problem at the initial stage of motion,'' Doklady Physics, {\bf60} (8), 368--371 (2015).

\bibitem{Morz}% Н. Ф. Морозов, А. К. Беляев, П. Е. Товстик, Т. П. Товстик,  ДАН,  {\bf 465}, 302 (2015). 
N. F. Morozov, A. K. Belyaev, P. E. Tovstik, T. P. Tovstik, ``
Initial stage of motion in the Lavrent’ev–Ishlinskii problem on longitudinal shock on a rod,'' Doklady Physics, {\bf60} (11|, 519–523 (2015).

\bibitem{Pai} M. A. Il'gamov, ``The interaction of instabilities in a hydroelastic system,'' // J. Appl. Math. and  Mech. {\bf 60}, p. 400--408 (2016).

\bibitem{Tih}% A. Н. Тихонов, В. Я. Арсенин, {\it Методы решения некорректных задач} (Наука , 1974).
A. N. Tikhonov, V. Y. Arsenin, 
{\it Methods for solving ill-posed problems} (John Wiley and Sons, 1977).

\bibitem{Lavm} %М. М. Лаврентьев, 
M. M. Lavrent'ev, ``On the Cauchy problem for Laplace equation,'' 
Izv. Akad. Nauk SSSR. Ser. Mat., 1956,{\bf20} (6), 819--842 (1956) [in russian].

\bibitem{Newd} D. J. Newman, ``Numerical method for solution of an elliptic Cauchy problem,''
 J. Math. and Phys., {\bf39}, 72--75 (1960/1961).

\bibitem{Iva}% B. K. Иванов,    
%Дифф. уравн. {\bf1}, 131 (1965).
 V. K. Ivanov, ``The Cauchy problem for the Laplace equation in an infinite strip // J. diff.
equation,'' {\bf 1}. 98–102 (1965).
 
\bibitem{Mas} %B. П. Маслов, 
%УМН {\bf XIII}, 183 (1968).
V. P. Maslov, `` The existence of a solution to an ill-posed problem is equivalent to the convergence of a regularization process,'' Uspechi Mat. Nauk {\bf23} (3), 183–184 (1980)[in russian].

\bibitem{Lali} %Р. Латтес, Ж. Л. Лионс, {\it Метод квазиобращения и его приложения} (Наука, 1974).
R. Lattes, J.-L. Lions, {\it The Method of Quasi-Reversibility: Applications to Partial Differential Equations}, (American Elsevier, 1969).

\bibitem{Meln} %И. В. Мельникова, У. А. Ануфриева, Совр. мат. Фунд. направления   
%{\bf 14}, 3 (2005).
I. V. Melnikova, U. A. Anufrieva, 
``Peculiarities and regularization of ill-posed Cauchy problems with differential operators,'' Journal of Mathematical Sciences
 {\bf 148} (4), 481–632 (2008).

\bibitem{Leo} A. C. Леонов, {\it Решение некорректно поставленных начальных задач: Очерк теории, практические алгоритмы и демонстрация в МАТЛАБ} (Либриком, 2009).

\bibitem{Tarh} N. N. Tarhanov, {\it The Cauchy problem for solutions of elliptic equations (Mathematical Topics} vol.7{\it)}
(Akademiс Verlag, 1995). 

\bibitem{Gursh}% А. В. Гуревич, А. Б. Шварцбург, ЖЭТФ {\bf 58}, 2012 (1970).
%Sov. Phys. JETP {\bf 31}, 1084 (1970)
 A.V. Gurevich, A.B. Shvartsburg,
``Exact Solutions of the Equations of Nonlinear Geometric Optics,'' Sov. Phys. JETP {\bf 31} (6), 1084--1089 (1970).


\bibitem{Kudb} %B. Р. Кудашев, Б. И. Сулейманов, %Особенности некоторых типичных процессов самопроизвольного
%падения интенсивности в неустойчивых средах.
%\type\published \lookin
%Письма в ЖЭТФ, {\bf 62}, 358, 1995.
V.R. Kudashev, B.I. Suleimanov,
``Characteristic features of some typical spontaneous intensivity collapse processes in unstable media,'' JETP Lett. {\bf 62} (4), 382--388 (1995).

\bibitem{Dub} B. Dubrovin, T. Grava, C. Klein, 
``On universality of critical behaviour in the critical behaviour in the focusing nonlinear Schr\"odinger equation, elliptic umbilic catstrophe and the tritonque to the Painlev\'e-I equation,''. J. Nonlinear Sc. {\bf 19} (1), 57--94 (2009).



\bibitem{Bkud}% B. Р. Кудашев, Б. И. Сулейманов,%Малоамплитудные дисперсионные колебания на фоне приближения
%нелинейной геометрической оптики.
%\type\published
%\lookin \lookin
%ТМФ, {\bf118}, 413, 1999.
 V.R. Kudashev, B.I. Suleimanov, 
``Small-amplitude dispersion oscillations on the background of the nonlinear geometric optics approximation,'' Theoret. and Math. Phys. {\bf118} (3), 325--332 (1999).

\bibitem{Tovb} M. Bertola, A. Tovbis, Coom. Pure and Appl. Math. {\bf 66}, 678 (2013).
``Universality for the focusing nonlinear Schr\"odinger equation at the gradient catastrophe point: rational breathers and poles of the tritronquйe solution to Painlev\'e I,'' Comm. Pure Appl. Math. {\bf 66} (5), 678--752 (2013). 

\bibitem{Gus}% В. И. Арнольд, А. Н. Варченко, С. М. Гусейн-Заде, {\it Особенности дифференцируемых отображений. Классификация критических точек, каустик и волновых фронтов}, (Наука, 1982).
V. I. Arnol'd, A. N. Varchenko, S. M. Gusein-Zade, {Singularities of Differentiable Mappings} V.1 
(Birkha\"user, Berlin, 1985).

\bibitem{Hardy} G. H. Hardy, ``On the asymptotic value of  certain integrals,'' Mess. Math. {\bf 46}, 70--73 (1917).

\bibitem{Rie} Э Я. Риекстыньш, {\it Асимптотические разложения интегралов} (Зинатне, 1977).

\bibitem{Ilam}% А. М. Ильин, {\it Согласование асимптотических разложений решений краевых задач} (Наука, 1989).
 A. M. Il'in, {\it Matching of asymptotic expansions of solutions of boundary value problems} Translations of Mathemathiclal Monographs, V. 102 (AMS, 1992).


\bibitem{Bab}% В. М. Бабич, В. С. Булдырев, {\it Асимптотические методы в задачах теории дифракции коротких волн. Метод эталонных задач} (Наука, 1972).
V. M. Babich, V. S. Buldyrev{\it Asymptotic methods in shortwave diffraction theory} (Springer, 1991)

\bibitem{Bryp} Ю. А. Брычков, А. П. Прудников {\it Интегральные преобразования обобщенных функций} (Наука, 1977).

%\bibitem{Erd} А. Эрдейи, {\it Асимптотические разложения}, ФИЗМАТЛИТ, Москва (1962). с. 1.

%\bibitem{Bab} Ю. А. Брычков, А. П. Прудников {\it Интегральные преобразования обобщенных функций}, Наука, Mocква (1977), с. 1.

\bibitem{Back} N. G. Backhoom, ``Asymptotic expansion of the function $F_k(x)=\int\limits_{0}^{\infty}e^{-u^k+xu}du$,'' Proc. London Math.Soc. {\bf 33} (2), 83--100 (1933).

\bibitem{Tolst} А. П. Прудников, Ю. А. Брычков, О. И. Маричев, {\it Интегралы и ряды} (Наука, 1981).

\bibitem{Beit} %Г. Бейтмен, А. Эрдейи {\it Таблицы интегральных преобразований}, 
%Том 1 (Наука,1969).
 H. Bateman, A. Erd\'elyi and other,  
{\it Tables of integral transforms} V.1 (McGraw--Hill, 1954)

\end{thebibliography}
\end{document}